Phase Factors in Singular Value Decomposition and Schmidt Decomposition

Chu Ryang Wie
University at Buffalo, State University of New York,
Department of Electrical Engineering, 230B Davis Hall, Buffalo, New York, USA

ABSTRACT
In singular value decomposition (*SVD*) of a complex matrix *A*, the singular vectors or the eigenvectors of $AA^\dagger$ and $A^\dagger A$ are unique up to complex phase factors. Thus, the two unitary matrices in *SVD* are unique up to diagonal matrices of phase factors, the phase-factor matrices. Also, the product of these two phase-factor matrices, or the product of phase factors of the corresponding singular vectors with the same singular value, is unique. In the Schmidt decomposition, a phase-factor matrix is a phase rotation operator acting on a subsystem alone. We summarize here three simple steps to consistently carry out the *SVD* and the Schmidt decomposition including the phase factors.

INTRODUCTION

Singular value decomposition (*SVD*) is widely used in numerical linear algebra. *SVD* is also used in obtaining the Schmidt decomposition of a bipartite quantum state in the quantum computing [1]. In a quantum computing course, the Schmidt decomposition of a bipartite quantum state is a necessary topic before the state purification is discussed in which a mixed state of a subsystem is converted to a pure state of a composite system in which the reduced density matrix for the subsystem remains the same.

The *SVD* theorem [2] may be stated as: "*Every m×n complex matrix A has an SVD, $A=UDV^\dagger$ where U is an m×m unitary matrix, V is an n×n unitary matrix, and D is a diagonal matrix with real nonnegative entries $\sigma_j$, called the singular values, arranged in nonincreasing order. Furthermore, the singular values $\{\sigma_j\}$ are uniquely determined, and if A is square and the $\sigma_j$ are distinct, the left and right singular vectors $|u_j\rangle$ and $|v_j\rangle$ are uniquely determined up to complex scalar factors of absolute value 1.*" Here, the singular vectors $|u_j\rangle$ and $|v_j\rangle$ are the eigenvectors of $AA^\dagger$ and $A^\dagger A$, respectively, corresponding to the eigenvalue $\sigma_j^2$.

When applied routinely, the *SVD* can lead to the same results for different matrices, such as for real symmetric, or complex *Hermitian*, matrices *A* and *B* where $A^2=B^2$. In this case, the *SVD* can yield the same *U*, *D,* and *V*. But, $UDV^\dagger$ may be equal to one matrix, not the other. Such examples are

$$A = \begin{pmatrix} 2 & 1 \\ 1 & 2 \end{pmatrix},\ B = \begin{pmatrix} 1 & 2 \\ 2 & 1 \end{pmatrix} \xRightarrow{SVD} U = \frac{1}{\sqrt{2}}\begin{pmatrix} 1 & 1 \\ 1 & -1 \end{pmatrix} = V,\ D = \begin{pmatrix} 3 & 0 \\ 0 & 1 \end{pmatrix}$$

or,

$$A = \begin{pmatrix} 3 & -2i \\ 2i & 3 \end{pmatrix},\ B = \begin{pmatrix} 2 & -3i \\ 3i & 2 \end{pmatrix} \xRightarrow{SVD} U = \frac{1}{\sqrt{2}}\begin{pmatrix} 1 & 1 \\ i & -i \end{pmatrix} = V,\ D = \begin{pmatrix} 5 & 0 \\ 0 & 1 \end{pmatrix}$$

In both cases, $UDV^\dagger$ is equal to *A*, but not to *B*. Students in a quantum computing class can encounter this kind of problems as assignments, spending a long fruitless hours on the unlucky case of matrix *B*. This result happens because the unitary matrices *U* and *V* are unique up to diagonal matrices of phase factors as a result of the fact that the eigenvectors are unique up to phase factors. The phase factors need be determined for each given matrix to do the *SVD* properly. One could be spending a significant amount of time before figuring out a solution. This note is written to present a streamlined and straightforward *SVD* procedure to save such a hassle and to lead to consistent results.

PROCEDURE FOR PHASE FACTORS IN SINGULAR VALUE DECOMPOSITION

Because the eigenvectors, *i.e.,* the singular vectors, are unique up to phase factors (*i.e.*, the complex scalar factors of modulus 1,) the unitary matrices in the *SVD* are unique up to diagonal matrices of the phase factors. The phase factor matrices must be determined for each matrix in the problem. In this section, we present a streamlined approach to determine the unitary matrices $U$ and $V$.

First, we relax the requirement of real nonnegative entries in the diagonal matrix, and allow the diagonal matrix to be complex. We determine the unitary matrices from any set of eigenvectors, disregarding the phase factors. Then, with these unitary matrices, we determine the complex entries for the complex diagonal matrix for the matrix in the problem. Next, the complex diagonal matrix is factorized into the real nonnegative diagonal matrix $D$ and a complex diagonal matrix of phase factors. The diagonal entries in the real diagonal matrix $D$ are the singular values. The diagonal matrix of phase factors is factored into two diagonal matrices, which are combined with the unitary matrices obtained earlier, to produce the final unitary matrices $U$ and $V$.

Therefore, the *SVD* procedure may be streamlined as follows: Let $A = U_0 D_0 V_0^\dagger = UDV^\dagger$ where $D_0$ is a complex diagonal matrix with complex entries $d_j$, and $D$ is a real diagonal matrix with nonnegative entries $\sigma_j = |d_j|$, the singular values of $A$. (*Step-1*) Disregard the phase factors and construct $U_0$ and $V_0$ from the eigenvectors of $AA^\dagger$ and $A^\dagger A$, respectively. (*Step-2*) Find $D_0$ by solving $U_0 D_0 V_0^\dagger = A$, and from the diagonal entries in $D_0$, find the entries (the absolute values) in $D$ and the phase factors. Factor the phase factor matrix into two phase factor matrices. (*Step-3*) Combine one phase factor matrix with $U_0$ to produce $U$, and the other phase factor matrix with $V_0$ for $V$. These three steps are detailed below:

*Step-1* Construct $U_0$ and $V_0$ from the eigenvectors $|u_j\rangle$ and $|v_j\rangle$ of $AA^\dagger$ and $A^\dagger A$, respectively, in order of the nonincreasing eigenvalues $\sigma_j^2$. Here, the phase factors are disregarded.
$$U_0 = (|u_0\rangle, |u_1\rangle, ..), \quad V_0 = (|v_0\rangle, |v_1\rangle, ..).$$

*Step-2* Find the complex diagonal matrix $D_0$, with $d_j$ as the $j^{th}$ entry, by solving $U_0 D_0 V_0^\dagger = \Sigma_j d_j |u_j\rangle\langle v_j| = A$. Let $d_j = \sigma_j \exp(i(\alpha_j + \beta_j))$ where the phase factors $\exp(i\alpha_j)$ and $\exp(i\beta_j)$ are set arbitrarily while their product is fixed. The modulus $|d_j| = \sigma_j$ is the $j^{th}$ singular value of $A$, and is the $j^{th}$ entry in the real diagonal matrix $D$.

$$D_0 = \begin{pmatrix} d_0 & 0 & 0 \\ 0 & d_1 & 0 \\ 0 & 0 & . \end{pmatrix} = \sum_j d_j |j\rangle\langle j| = \sum_j \sigma_j e^{i(\alpha_j+\beta_j)} |j\rangle\langle j| = D_a D D_b^\dagger,$$

$$D_a \equiv \begin{pmatrix} e^{i\alpha_0} & 0 & 0 \\ 0 & e^{i\alpha_1} & 0 \\ 0 & 0 & . \end{pmatrix} = \sum_j e^{i\alpha_j} |j\rangle\langle j|,$$

$$D_b \equiv \begin{pmatrix} e^{-i\beta_0} & 0 & 0 \\ 0 & e^{-i\beta_1} & 0 \\ 0 & 0 & . \end{pmatrix} = \sum_j e^{-i\beta_j} |j\rangle\langle j|,$$

$$D = \begin{pmatrix} |d_0| & 0 & 0 \\ 0 & |d_1| & 0 \\ 0 & 0 & . \end{pmatrix} = \begin{pmatrix} \sigma_0 & 0 & 0 \\ 0 & \sigma_1 & 0 \\ 0 & 0 & . \end{pmatrix} = \sum_j \sigma_j |j\rangle\langle j|$$

*Step-3* Finally, the phase factor matrices $D_a$ and $D_b$ with diagonal entries of $\exp(i\alpha_j)$ and $\exp(-i\beta_j)$, respectively, are combined with the unitary matrices $U_0$ and $V_0$ to produce $U=U_0 D_a$ and $V=V_0 D_b$. That is, the $j^{th}$ column vectors are $exp(i\alpha_j)|u_j\rangle$ and $exp(-i\beta_j)|v_j\rangle$ for $U$ and $V$, respectively. The SVD is $A = UDV^\dagger = U_0D_aDD_b^\dagger V_0^\dagger = \Sigma_j\sigma_j exp(i\alpha_j)|u_j\rangle\langle v_j|exp(i\beta_j)$.

For each nonzero entry $d_j$ in $D_0$, define:
$$e^{i(\alpha_j+\beta_j)} \equiv \frac{d_j}{|d_j|} = \frac{d_j}{\sigma_j}$$

The SVD is complete once $U=U_0D_a$ and $V=V_0D_b$ are obtained as:

a. $U = \left(e^{i\alpha_0}|u_0\rangle,\ e^{i\alpha_1}|u_1\rangle,..\right),\ V = \left(e^{-i\beta_0}|v_0\rangle,\ e^{-i\beta_1}|v_1\rangle,..\right)$

All phase factors may be assigned to one matrix as below:

b. $U = \left(\frac{d_0}{|d_0|}|u_0\rangle, \frac{d_1}{|d_1|}|u_1\rangle,..\right),\ V = (|v_0\rangle, |v_1\rangle,..)$

c. $U = (|u_0\rangle, |u_1\rangle,..),\ V = \left(\frac{d_0^*}{|d_0|}|v_0\rangle, \frac{d_1^*}{|d_1|}|v_1\rangle,..\right)$

Where, $d_j^*$ is a complex conjugate of $d_j$.

SCHMIDT DECOMPOSITION

The Schmidt decomposition for a bipartite state is as follows: For a bipartite quantum state $|\psi\rangle$, let
$$|\psi\rangle = \sum_{j,k} a_{jk}\ |j\rangle|k\rangle = \sum_i \sigma_i|i_A\rangle|i_B\rangle$$

Where, $|j\rangle$ and $|k\rangle$ are the computational basis states, and $|i_A\rangle$ and $|i_B\rangle$ are the orthonormal Schmidt basis states for systems $A$ and $B$, respectively. The singular value $\sigma_i$ of matrix $(a_{jk})$ is the $i^{th}$ Schmidt coefficient. For a two-qubit state, the Schmidt basis can be obtained as follows. With the SVD of the matrix $(a_{jk}) = UDV^\dagger$ as obtained in the three steps of SVD above, we find the Schmidt bases:
$$\begin{pmatrix}|0_A\rangle\\|1_A\rangle\end{pmatrix} = U^T\begin{pmatrix}|0\rangle\\|1\rangle\end{pmatrix},$$
$$\begin{pmatrix}|0_B\rangle\\|1_B\rangle\end{pmatrix} = V^\dagger\begin{pmatrix}|0\rangle\\|1\rangle\end{pmatrix}.$$

Here, $U^T$ is the transpose of $U$. The *ket*-vector entries in the 2-dimensional column vectors are treated each as a single entity (or entry). The matrix products are carried out on the right hand side.

The phase factor matrices, $D_a$ and $D_b$, act on each qubit alone as a phase rotation operator, a unitary gate.
$$\begin{pmatrix}|0_A\rangle\\|1_A\rangle\end{pmatrix} = U^T\begin{pmatrix}|0\rangle\\|1\rangle\end{pmatrix} = D_a^T U_0^T\begin{pmatrix}|0\rangle\\|1\rangle\end{pmatrix},$$
$$\begin{pmatrix}|0_B\rangle\\|1_B\rangle\end{pmatrix} = V^\dagger\begin{pmatrix}|0\rangle\\|1\rangle\end{pmatrix} = D_b^\dagger V_0^\dagger\begin{pmatrix}|0\rangle\\|1\rangle\end{pmatrix}.$$

The different constructions of $U=U_0D_a$ and $V=V_0D_b$ from the SVD *step-3* above, will lead to different Schmidt basis states $|i_A\rangle$ and $|i_B\rangle$, different by the phase factors, but they combine and produce an equally valid Schmidt decomposition $|\psi\rangle = \Sigma_j\ \sigma_j\ |j_A\rangle|j_B\rangle$.

CALCULATION EXAMPLES

Here, we discuss several SVD examples and a Schmidt decomposition example. For $2\times 2$ matrices, where convenient, we shall use in the calculations the Pauli matrices $X$, $Y$, and $Z$ and their eigenvectors. Pauli matrices are *Hermitian* and self-inverse.

(1) Find the SVD for $B = \begin{pmatrix}1 & 2\\2 & 1\end{pmatrix} = I + 2X$

This matrix B is one of the aforementioned examples for which a routine procedure of SVD can lead to some confusing results at the end. Here, $I$ is a $2 \times 2$ identity matrix, and $X$ is a Pauli matrix. Using the *Hermitian* and self-inverse properties of Pauli matrices where $X^2 = I$, we get $BB^\dagger = B^\dagger B = B^2 = (I+2X)^2 = 5I + 4X$. Pauli $X$ has eigenvectors $|+\rangle$ and $|-\rangle$ and eigenvalues +1 and -1, respectively, with the following matrix representations:

$$|+\rangle = \frac{1}{\sqrt{2}}\begin{pmatrix}1\\1\end{pmatrix}, \quad |-\rangle = \frac{1}{\sqrt{2}}\begin{pmatrix}1\\-1\end{pmatrix}, \quad |0\rangle = \begin{pmatrix}1\\0\end{pmatrix}, \quad |1\rangle = \begin{pmatrix}0\\1\end{pmatrix}$$

Therefore, for $B^2$, the eigenvalues are 5+4=9 and 5-4=1, with eigenvectors $|+\rangle$ and $|-\rangle$, respectively. Disregarding the phase factors, let $U_0 = (|+\rangle, |-\rangle) = V_0$,

$$U_0 D_0 V_0^\dagger = (d_0|+\rangle, d_1|-\rangle)\begin{pmatrix}\langle+|\\\langle-|\end{pmatrix} = d_0|+\rangle\langle+| + d_1|-\rangle\langle-| = |d_0|\frac{d_0}{|d_0|}|+\rangle\langle+| + |d_1|\frac{d_1}{|d_1|}|-\rangle\langle-|$$

Using the identities: $|+\rangle\langle+| = \frac{1}{2}(I+X)$ and $|-\rangle\langle-| = \frac{1}{2}(I-X)$, and solving for $d_0$ and $d_1$,

$$U_0 D_0 V_0^\dagger = \frac{d_0+d_1}{2}I + \frac{d_0-d_1}{2}X = I + 2X$$

We get $d_0 = 3$, $d_1 = -1$, and the phase factors are $d_0/|d_0| = \exp(i(\alpha_0+\beta_0)) = 1$ and $d_1/|d_1| = \exp(i(\alpha_1+\beta_1)) = -1$. The real diagonal matrix $D$ with nonnegative (singular value) entries is $\begin{pmatrix}3 & 0\\0 & 1\end{pmatrix}$, and we could use for SVD any of the following $U$ and $V$:

a. $U = (e^{i\alpha_0}|+\rangle, e^{i\alpha_1}|-\rangle) = \frac{1}{\sqrt{2}}\begin{pmatrix}1 & i\\1 & -i\end{pmatrix}$, $V = (e^{-i\beta_0}|+\rangle, e^{-i\beta_1}|-\rangle) = \frac{1}{\sqrt{2}}\begin{pmatrix}1 & -i\\1 & i\end{pmatrix}$,
   $\alpha_0 = 0 = \beta_0$, and $\alpha_1 = \pi/2 = \beta_1$.

b. $U = \left(\frac{d_0}{|d_0|}|+\rangle, \frac{d_1}{|d_1|}|-\rangle\right) = \frac{1}{\sqrt{2}}\begin{pmatrix}1 & -1\\1 & 1\end{pmatrix}$, $V = (|+\rangle, |-\rangle) = \frac{1}{\sqrt{2}}\begin{pmatrix}1 & 1\\1 & -1\end{pmatrix}$

c. $U = (|+\rangle, |-\rangle) = \frac{1}{\sqrt{2}}\begin{pmatrix}1 & 1\\1 & -1\end{pmatrix}$, $V = \left(\frac{d_0^*}{|d_0|}|+\rangle, \frac{d_1^*}{|d_1|}|-\rangle\right) = \frac{1}{\sqrt{2}}\begin{pmatrix}1 & -1\\1 & 1\end{pmatrix}$

(2) Find the SVD for $B = \begin{pmatrix}2 & -3i\\3i & 2\end{pmatrix} = 2I + 3Y$

This matrix $B$ was the second of the aforementioned examples for which a routine SVD can lead to a confusing result at the end. Here, $Y$ is a Pauli matrix. $BB^\dagger = B^\dagger B = B^2 = (2I+3Y)^2 = 13I + 12Y$. The eigenvectors for $Y$ are $|+\rangle_y$ and $|-\rangle_y$ with eigenvalues +1 and -1, respectively, with the following matrix representations:

$$|+\rangle_y = \frac{1}{\sqrt{2}}\begin{pmatrix}1\\i\end{pmatrix}, \quad |-\rangle_y = \frac{1}{\sqrt{2}}\begin{pmatrix}1\\-i\end{pmatrix}.$$

Therefore, for $B^2$, the eigenvalues are 13+12=25 and 13-12=1, and the eigenvectors are the same as for $Y$, $|+\rangle_y$ and $|-\rangle_y$. Disregarding the phase factors, set $U_0 = (|+\rangle_y, |-\rangle_y) = V_0$, and determine $D_0$ by solving $U_0 D_0 V_0^\dagger = B$.

$$U_0 D_0 V_0^\dagger = (d_0|+\rangle_y, d_1|-\rangle_y)\begin{pmatrix}\langle+|_y\\\langle-|_y\end{pmatrix} = d_0|+\rangle_y\langle+| + d_1|-\rangle_y\langle-|$$

$$= |d_0|\frac{d_0}{|d_0|}|+\rangle_y\langle+| + |d_1|\frac{d_1}{|d_1|}|-\rangle_y\langle-|$$

Using the identities: $|+\rangle_y\langle+| = \frac{1}{2}(I+Y)$ and $|-\rangle_y\langle-| = \frac{1}{2}(I-Y)$, solve the following,

$$U_0 D_0 V_0^\dagger = \frac{d_0+d_1}{2}I + \frac{d_0-d_1}{2}Y = 2I + 3Y$$

We get $d_0 = 5$, $d_1 = -1$, $d_0/|d_0|=\exp(i(\alpha_0+\beta_0))=1$, and $d_1/|d_1|=\exp(i(\alpha_1+\beta_1))=-1$. The real diagonal matrix $D$ is $D = \begin{pmatrix}5 & 0\\0 & 1\end{pmatrix}$. We could get any of the following set of $U$ and $V$:

a. $U = (e^{i\alpha_0}|+\rangle_y, e^{i\alpha_1}|-\rangle_y) = \frac{1}{\sqrt{2}}\begin{pmatrix}1 & i\\i & 1\end{pmatrix}$,
   $V = (e^{-i\beta_0}|+\rangle_y, e^{-i\beta_1}|-\rangle_y) = \frac{1}{\sqrt{2}}\begin{pmatrix}1 & -i\\i & -1\end{pmatrix}$, $\alpha_0 = 0 = \beta_0$, and $\alpha_1 = \pi/2 = \beta_1$.

b. $U = \left(\frac{d_0}{|d_0|}|+\rangle_y, \frac{d_1}{|d_1|}|-\rangle_y\right) = \frac{1}{\sqrt{2}}\begin{pmatrix} 1 & -1 \\ i & i \end{pmatrix}$, $V = (|+\rangle_y, |-\rangle_y) = \frac{1}{\sqrt{2}}\begin{pmatrix} 1 & 1 \\ i & -i \end{pmatrix}$

c. $U = (|+\rangle_y, |-\rangle_y) = \frac{1}{\sqrt{2}}\begin{pmatrix} 1 & 1 \\ i & -i \end{pmatrix}$, $V = \left(\frac{d_0^*}{|d_0|}|+\rangle_y, \frac{d_1^*}{|d_1|}|-\rangle_y\right) = \frac{1}{\sqrt{2}}\begin{pmatrix} 1 & -1 \\ i & i \end{pmatrix}$

(3) Find the *SVD* for $A = \begin{pmatrix} 1 & 1 & i \\ 1 & -1 & i \end{pmatrix}$

The same *SVD* procedure is applied to this 2×3 complex rectangular matrix *A*. First, we factorize *A* into $U_0 D_0 V_0^\dagger$ for the unitary matrices $U_0$ and $V_0$ formed from the eigenvectors of $AA^\dagger$ and $A^\dagger A$, respectively, and the complex diagonal matrix $D_0$ with the diagonal entries $d_0$, $d_1$, etc. determined by solving $U_0 D_0 V_0^\dagger = A$. $AA^\dagger$, a 2×2 matrix, has eigenvalues 4 and 2, with eigenvectors $|u_0\rangle=|+\rangle$ and $|u_1\rangle=|-\rangle$, respectively. The 3×3 matrix $A^\dagger A$ has eigenvalues 4, 2 and 0, with the following eigenvectors, respectively.

$$|v_0\rangle = \frac{1}{\sqrt{2}}\begin{pmatrix} i \\ 0 \\ 1 \end{pmatrix}, \quad |v_1\rangle = \begin{pmatrix} 0 \\ 1 \\ 0 \end{pmatrix}, \quad |v_2\rangle = \frac{1}{\sqrt{2}}\begin{pmatrix} -i \\ 0 \\ 1 \end{pmatrix}$$

We then find $d_0$ and $d_1$ by solving $U_0 D_0 V_0^\dagger = A$, and $d_2$ is zero.

$$U_0 D_0 V_0^\dagger = d_0|u_0\rangle\langle v_0| + d_1|u_1\rangle\langle v_1| = \begin{pmatrix} -i d_0/2 & d_1/\sqrt{2} & d_0/2 \\ -i d_0/2 & -d_1/\sqrt{2} & d_0/2 \end{pmatrix} = \begin{pmatrix} 1 & 1 & i \\ 1 & -1 & i \end{pmatrix}$$

This gives $d_0 = 2i$ and $d_1=\sqrt{2}$. In step-2 of *SVD*, the phase factors are: $d_0/|d_0|=i$, and $d_1/|d_1|=1$ which may be absorbed into either *U* or *V*, or may be split between the two. Here, the phase factors are absorbed in *U*. Hence, $U = (i|u_0\rangle, |u_1\rangle)$, $V=(|v_0\rangle, |v_1\rangle, |v_2\rangle)$ and the 2×3 diagonal matrix *D* has real nonnegative entries 2 and $\sqrt{2}$ in the descending order. Thus, $UDV^\dagger$ is an *SVD* of this rectangular matrix *A*.

(4) Perform the Schmidt decomposition for $|\psi\rangle = \frac{2i |00\rangle + |01\rangle + |10\rangle + 2 |11\rangle}{\sqrt{10}}$

The matrix to do the *SVD* on is $A = \frac{1}{\sqrt{10}}\begin{pmatrix} 2i & 1 \\ 1 & 2 \end{pmatrix}$. Matrices $AA^\dagger$ and $A^\dagger A$ have the same eigenvalues $|d_0|^2 = (5+2\sqrt{2})/10$ and $|d_1|^2 = (5-2\sqrt{2})/10$, and the following eigenvectors, respectively.

$$|u_0\rangle = \frac{1}{2}\begin{pmatrix} \sqrt{2} \\ 1-i \end{pmatrix}, |u_1\rangle = \frac{1}{2}\begin{pmatrix} -\sqrt{2} \\ 1-i \end{pmatrix}, |v_0\rangle = \frac{1}{2}\begin{pmatrix} \sqrt{2} \\ 1+i \end{pmatrix}, |v_1\rangle = \frac{1}{2}\begin{pmatrix} -\sqrt{2} \\ 1+i \end{pmatrix}$$

The diagonal entries $d_0$ and $d_1$ to the complex diagonal matrix $D_0$ are found by solving $d_0|u_0\rangle\langle v_0|+d_1|u_1\rangle\langle v_1| = A$.

$$d_0 = \frac{\sqrt{2} + (4+\sqrt{2})i}{2\sqrt{10}}, \quad d_1 = \frac{-\sqrt{2} + (4-\sqrt{2})i}{2\sqrt{10}}, \quad e^{i(\alpha_j+\beta_j)} = \frac{d_j}{|d_j|}, \quad j = 0, 1$$

They yield the phase factors $exp(i(\alpha_j+\beta_j))$ for $j=0$ and $1$ where the division between $\alpha_j$ and $\beta_j$ is arbitrary for each *j*. For $A = UDV^\dagger$ with the diagonal matrix *D* with real nonnegative entries $|d_0|$ and $|d_1|$, the unitary matrices *U* and *V* are

$$U = (e^{i\alpha_0}|u_0\rangle, e^{i\alpha_1}|u_1\rangle) = \frac{1}{2}\begin{pmatrix} \sqrt{2}e^{i\alpha_0} & -\sqrt{2}e^{i\alpha_1} \\ (1-i)e^{i\alpha_0} & (1-i)e^{i\alpha_1} \end{pmatrix},$$

$$V = (e^{-i\beta_0}|v_0\rangle, e^{-i\beta_1}|v_1\rangle) = \frac{1}{2}\begin{pmatrix} \sqrt{2}e^{-i\beta_0} & -\sqrt{2}e^{-i\beta_1} \\ (1+i)e^{-i\beta_0} & (1+i)e^{-i\beta_1} \end{pmatrix}$$

The Schmidt bases are found in terms of the computational bases as

$$\begin{pmatrix} |0_A\rangle \\ |1_A\rangle \end{pmatrix} = U^T \begin{pmatrix} |0\rangle \\ |1\rangle \end{pmatrix} = \frac{1}{2}\begin{pmatrix} \sqrt{2}e^{i\alpha_0}|0\rangle + (1-i)e^{i\alpha_0}|1\rangle \\ -\sqrt{2}e^{i\alpha_1}|0\rangle + (1-i)e^{i\alpha_1}|1\rangle \end{pmatrix} = \begin{pmatrix} e^{i\alpha_0} & 0 \\ 0 & e^{i\alpha_1} \end{pmatrix}\frac{1}{2}\begin{pmatrix} \sqrt{2}|0\rangle + (1-i)|1\rangle \\ -\sqrt{2}|0\rangle + (1-i)|1\rangle \end{pmatrix},$$

$$\begin{pmatrix} |0_B\rangle \\ |1_B\rangle \end{pmatrix} = V^\dagger \begin{pmatrix} |0\rangle \\ |1\rangle \end{pmatrix} = \frac{1}{2}\begin{pmatrix} \sqrt{2}e^{i\beta_0}|0\rangle + (1-i)e^{i\beta_0}|1\rangle \\ -\sqrt{2}e^{i\beta_1}|0\rangle + (1-i)e^{i\beta_1}|1\rangle \end{pmatrix} = \begin{pmatrix} e^{i\beta_0} & 0 \\ 0 & e^{i\beta_1} \end{pmatrix}\frac{1}{2}\begin{pmatrix} \sqrt{2}|0\rangle + (1-i)|1\rangle \\ -\sqrt{2}|0\rangle + (1-i)|1\rangle \end{pmatrix},$$

Here, each phase factor matrix is a unitary phase-rotation operator acting on the subsystem alone. This is notable because if $\Sigma_i\sigma_i|i_A\rangle|i_B\rangle$ is the Schmidt decomposition for a state $|\psi\rangle$, then $\Sigma_i\sigma_i(U|i_A\rangle)|i_B\rangle$ is the Schmidt decomposition for $U|\psi\rangle$, where $U$ is a unitary operator acting on system $A$ alone [1]. In fact, the unitary operators $U_0$ and $V_0$, or $U$ and $V$, from the *SVD* are unitary gates acting on each qubit alone. The Schmidt decomposition with the phase factors arbitrarily split between the two qubits works as below.

$$
\begin{aligned}
|\psi\rangle &= |d_0||0_A\rangle|0_B\rangle + |d_1||1_A\rangle|1_B\rangle \\
&= \frac{|d_0|}{4}e^{i(\alpha_0+\beta_0)}(\sqrt{2}|0\rangle + (1-i)|1\rangle)(\sqrt{2}|0\rangle + (1-i)|1\rangle) \\
&\quad + \frac{|d_1|}{4}e^{i(\alpha_1+\beta_1)}(-\sqrt{2}|0\rangle + (1-i)|1\rangle)(-\sqrt{2}|0\rangle + (1-i)|1\rangle) \\
&= \frac{d_0}{4}(\sqrt{2}|0\rangle + (1-i)|1\rangle)(\sqrt{2}|0\rangle + (1-i)|1\rangle) \\
&\quad + \frac{d_1}{4}(-\sqrt{2}|0\rangle + (1-i)|1\rangle)(-\sqrt{2}|0\rangle + (1-i)|1\rangle) \\
&= \frac{d_0+d_1}{2}(|00\rangle - i|11\rangle) + \frac{d_0-d_1}{4}\sqrt{2}(1-i)(|01\rangle+|10\rangle) \\
&= \frac{2i|00\rangle + |01\rangle + |10\rangle + 2|11\rangle}{\sqrt{10}}
\end{aligned}
$$

Again, the diagonal complex matrix of phase factors $d_0/|d_0|$ and $d_1/|d_1|$ can be merged either with $U_0$ to produce $U$ (while $V=V_0$), its complex conjugate merged with $V_0$ to produce $V$ (while $U=U_0$), or the phase factors can be arbitrarily split between the two unitary matrices $U$ and $V$. Here we followed the third way.

PROCEDURE SUMMARY

The *SVD* of a complex matrix $A$ was performed by following three simple steps to set the proper phase factors in the unitary matrices: $A = U_0 D_0 V_0^\dagger = UDV^\dagger$ where $D_0$ is a complex diagonal matrix, and $D$ is a real diagonal matrix where the entries are the singular values or the absolute values of entries in $D_0$. The column vectors of a unitary matrix are unique up to phase factors. Therefore, given a complex matrix $A$, the complex diagonal matrix $D_0$ is used to find the correct phase factors for the column vectors of $U$ and $V$.

*Step-1* Disregarding the phase factors, construct $U_0$ and $V_0$ from the eigenvectors $|u_j\rangle$ and $|v_j\rangle$ of $AA^\dagger$ and $A^\dagger A$, respectively, in order of the decreasing eigenvalue.
$$U_0 = (|u_0\rangle, |u_1\rangle, ..), \quad V_0 = (|v_0\rangle, |v_1\rangle, ..).$$

*Step-2* Find the entries $d_j$ of the complex diagonal matrix $D_0$ by solving $U_0 D_0 V_0^\dagger = \Sigma_j d_j|u_j\rangle\langle v_j| = A$. Find the phase factors $\exp(i(\alpha_j+\beta_j)) = d_j/|d_j| = d_j/\sigma_j$. Here $d_j$ is the $j^{th}$ entry of $D_0$ and $\sigma_j=|d_j|$ is the $j^{th}$ entry of $D$ and is a singular value. $D_0 = \Sigma_j d_j |j\rangle\langle j|$, and $D = \Sigma_j \sigma_j |j\rangle\langle j|$.

*Step-3* Finally, combine the phase factors with the column vectors of $U_0$ and $V_0$, and construct the unitary matrices $U$ and $V$, respectively.
   a. $U = (e^{i\alpha_0}|u_0\rangle, \ e^{i\alpha_1}|u_1\rangle, ..), \quad V = (e^{-i\beta_0}|v_0\rangle, \ e^{-i\beta_1}|v_1\rangle, ..)$

Note that the phase factors can be divided arbitrarily as long as their product remains the same. Therefore, a simple construction may be:

    b. $U = \left(\frac{d_0}{|d_0|}|u_0\rangle, \frac{d_1}{|d_1|}|u_1\rangle, ..\right), \ V = (|v_0\rangle, |v_1\rangle, ..)$

    c. $U = (|u_0\rangle, |u_1\rangle, ...), \ V = \left(\frac{{d_0}^*}{|d_0|}|v_0\rangle, \frac{{d_1}^*}{|d_1|}|v_1\rangle, ..\right)$

SUMMARY

A streamlined three-step procedure is reported here to set the phase factors correctly in the unitary matrices in the *SVD* and in the Schmidt decomposition. This procedure can produce consistent results for the SVD and Schmidt decomposition, including the phase factors, of an arbitrary complex matrix.

REFERENCES

1. Quantum Computation and Quantum Information, Michael A. Nielsen and Isaac L. Chuang, pp.78-79, pp.109-119, Cambridge University Press, 2010
2. Numerical Linear Algebra, Lloyd N. Trefethen and David Bau, III, pp.25-30, Society For Industrial And Applied Mathematics, 1997